\documentclass[english,oneside]{amsart}
\usepackage[T1]{fontenc}
\usepackage[latin1]{inputenc}
\usepackage{amssymb}

\makeatletter
 \theoremstyle{plain}
\newtheorem{thm}{Theorem}[section]
  \theoremstyle{plain}
  \newtheorem{lem}[thm]{Lemma}

\usepackage{babel}
\makeatother
\begin{document}

\title{Scaled Asymptotics For Some $q$-Series  }

\author{Ruiming Zhang}


\email{ruimingzhang@yahoo.com}

\begin{abstract}
In this work we investigate the asymptotics for Euler's $q$-Exponential
$E_{q}(z)$, $q$-Gamma function $\Gamma_{q}(z)$, Ramanujan's function
$A_{q}(z)$, Jackson's $q$-Bessel function $J_{\nu}^{(2)}$(z;q)
of second kind, Stieltjes-Wigert orthogonal polynomials $S_{n}(x;q)$
and $q$-Laguerre polynomials $L_{n}^{(\alpha)}(x;q)$ as $q$ approaching
$1$. 
\end{abstract}

\date{November 30, 2006}

\curraddr{School of Mathematics\\
Guangxi Normal University\\
Guilin City, Guangxi 541004\\
P. R. China.}

\address{Adjunct Professor \\
Binzhou Vocational College\\
533 Bohai 9 Road\\
Binzhou City, Shandong 256624\\
P. R. China.}

\maketitle

\section{Introduction}

Euler's $q$-Exponential $E_{q}(z)$, $q$-Gamma function $\Gamma_{q}(z)$,
Ramanujan's function $A_{q}(z)$, Jackson's $q$-Bessel function $J_{\nu}^{(2)}$(z;q)
of second kind, Stieltjes-Wigert orthogonal polynomials $S_{n}(x;q)$
and $q$-Laguerre polynomials $L_{n}^{(\alpha)}(x;q)$ are very important
examples in $q$-series, \cite{Andrews4,Gasper,Ismail2,Koekoek}.
They are used widely in other branches of mathematics and physics.
In this work we will present a method to derive asymptotic formulas
for these functions as $q$ approaching $1$ via Jacobi theta functions.

For any complex number $a$ and parameter $0<q<1$, we define \cite{Andrews4,Gasper,Ismail2,Koekoek}\begin{equation}
(a;q)_{\infty}:=\prod_{k=0}^{\infty}(1-aq^{k})\label{eq:1.1}\end{equation}
 and the $q$-shifted factorial of $a$ is defined by

\begin{equation}
(a;q)_{n}:=\frac{(a;q)_{\infty}}{(aq^{n};q)_{\infty}},\quad n\in\mathbb{Z}.\label{eq:1.2}\end{equation}
 We also use the following short-hand notation\begin{equation}
(a_{1},\dots,a_{m};q)_{n}=\prod_{k=1}^{m}(a_{k};q)_{n},\quad m\in\mathbb{N},\quad a_{1},\dots,a_{m}\in\mathbb{C},\quad n\in\mathbb{Z}.\label{eq:1.3}\end{equation}

\begin{lem}
\label{lem:1}Given any complex number $a$, assume that\begin{equation}
0<\frac{\left|a\right|q^{n}}{1-q}<\frac{1}{2}\label{eq:1.4}\end{equation}
 for some positive integer $n$. Then, \begin{equation}
\frac{(a;q)_{\infty}}{(a;q)_{n}}=(aq^{n};q)_{\infty}:=1+r_{1}(a;n)\label{eq:1.5}\end{equation}
 with\begin{equation}
\left|r_{1}(a;n)\right|\le\frac{2\left|a\right|q^{n}}{1-q}\label{eq:1.6}\end{equation}
 and \begin{equation}
\frac{(a;q)_{n}}{(a;q)_{\infty}}=\frac{1}{(aq^{n};q)_{\infty}}:=1+r_{2}(a;n)\label{eq:1.7}\end{equation}
 with \begin{equation}
\left|r_{2}(a;n)\right|\le\frac{2\left|a\right|q^{n}}{(1-q)}.\label{eq:1.8}\end{equation}

\end{lem}
\begin{proof}
From the $q$-binomial theorem \cite{Andrews4,Gasper,Ismail2,Koekoek}\begin{equation}
\frac{(az;q)_{\infty}}{(z;q)_{\infty}}=\sum_{k=0}^{\infty}\frac{(a;q)_{k}}{(q;q)_{k}}z^{k}\quad a,z\in\mathbb{C},\label{eq:1.9}\end{equation}
 and the inequality \begin{equation}
(q;q)_{k}\ge(1-q)^{k}\label{eq:1.10}\end{equation}
 for $k=0,1,...$, we obtain \begin{eqnarray}
r_{2}(a;n) & = & \sum_{k=0}^{\infty}\frac{\left(aq^{n}\right)^{k+1}}{(q;q)_{k+1}}\label{eq:1.11}\end{eqnarray}
 and \begin{align}
 & \left|r_{2}(a;n)\right|\le\sum_{k=0}^{\infty}\frac{\left(\left|a\right|q^{n}\right)^{k+1}}{(q;q)_{k+1}}\le\frac{\left|a\right|q^{n}}{(1-q)}\sum_{k=0}^{\infty}\left(\frac{\left|a\right|q^{n}}{1-q}\right)^{k}\le\frac{2\left|a\right|q^{n}}{(1-q)}.\label{eq:1.12}\end{align}
 Apply a limiting case of \eqref{eq:1.9}, \begin{equation}
(z;q)_{\infty}=\sum_{k=0}^{\infty}\frac{q^{k(k-1)/2}}{(q;q)_{k}}(-z)^{k}\quad z\in\mathbb{C},\label{eq:1.13}\end{equation}
 and the inequalities,\begin{equation}
\frac{1-q^{k}}{1-q}\ge kq^{k-1},\quad\frac{(q;q)_{k}}{(1-q)^{k}}\ge k!q^{k(k-1)/2},\quad\mbox{for }k=0,1,\dots\label{eq:1.14}\end{equation}
to obtain\begin{align}
r_{1}(a;n) & =\sum_{k=1}^{\infty}\frac{q^{k(k-1)/2}(-aq^{n})^{k}}{(q;q)_{k}},\label{eq:1.15}\end{align}
 and \begin{align}
 & \left|r_{1}(a;n)\right|\le\sum_{k=0}^{\infty}\frac{(\left|a\right|q^{n})^{k+1}}{(1-q)^{k+1}}\frac{(1-q)^{k+1}q^{k(k+1)/2}}{(q;q)_{k+1}}\label{eq:1.16}\\
 & \le\sum_{k=0}^{\infty}\frac{(\left|a\right|q^{n})^{k+1}}{(1-q)^{k+1}}\frac{1}{(k+1)!}\le\frac{\left|a\right|q^{n}}{1-q}\exp(1/2)<\frac{2\left|a\right|q^{n}}{1-q}.\nonumber \end{align}

\end{proof}
The Dedekind $\eta(\tau)$ is defined as \cite{Rademarcher}\begin{equation}
\eta(\tau):=e^{\pi i\tau/12}\prod_{k=1}^{\infty}(1-e^{2\pi ik\tau}),\label{eq:1.17}\end{equation}
 or\begin{equation}
\eta(\tau)=q^{1/12}(q^{2};q^{2})_{\infty},\quad q=e^{\pi i\tau},\quad\Im(\tau)>0.\label{eq:1.18}\end{equation}
 It has the transformation formula\begin{equation}
\eta\left(-\frac{1}{\tau}\right)=\sqrt{\frac{\tau}{i}}\eta(\tau).\label{eq:1.19}\end{equation}
 The Jacobi theta functions are defined as \begin{align}
\theta_{1}(z;q):=\theta_{1}(v|\tau) & :=-i\sum_{k=-\infty}^{\infty}(-1)^{k}q^{(k+1/2)^{2}}e^{(2k+1)\pi iv},\label{eq:1.20}\\
\theta_{2}(z;q):=\theta_{2}(v|\tau) & :=\sum_{k=-\infty}^{\infty}q^{(k+1/2)^{2}}e^{(2k+1)\pi iv},\label{eq:1.21}\\
\theta_{3}(z;q):=\theta_{3}(v|\tau) & :=\sum_{k=-\infty}^{\infty}q^{k^{2}}e^{2k\pi iv},\label{eq:1.22}\\
\theta_{4}(z;q):=\theta_{4}(v|\tau) & :=\sum_{k=-\infty}^{\infty}(-1)^{k}q^{k^{2}}e^{2k\pi iv},\label{eq:1.23}\end{align}
 where\begin{equation}
z=e^{2\pi iv},\quad q=e^{\pi i\tau},\quad\Im(\tau)>0.\label{eq:1.24}\end{equation}
 The Jacobi's triple product identities are\begin{align}
\theta_{1}(v|\tau) & =2q^{1/4}\sin\pi v(q^{2};q^{2})_{\infty}(q^{2}e^{2\pi iv};q^{2})_{\infty}(q^{2}e^{-2\pi iv};q^{2})_{\infty},\label{eq:1.25}\\
\theta_{2}(v|\tau) & =2q^{1/4}\cos\pi v(q^{2};q^{2})_{\infty}(-q^{2}e^{2\pi iv};q^{2})_{\infty}(-q^{2}e^{-2\pi iv};q^{2})_{\infty},\label{eq:1.26}\\
\theta_{3}(v|\tau) & =(q^{2};q^{2})_{\infty}(-qe^{2\pi iv};q^{2})_{\infty}(-qe^{-2\pi iv};q^{2})_{\infty},\label{eq:1.27}\\
\theta_{4}(v|\tau) & =(q^{2};q^{2})_{\infty}(qe^{2\pi iv};q^{2})_{\infty}(qe^{-2\pi iv};q^{2})_{\infty}.\label{eq:1.28}\end{align}
 The Jacobi $\theta$ functions satisfy transformations:\begin{align}
\theta_{1}\left(\frac{v}{\tau}\mid-\frac{1}{\tau}\right) & =-i\sqrt{\frac{\tau}{i}}e^{\pi iv^{2}/\tau}\theta_{1}\left(v\mid\tau\right),\label{eq:1.29}\\
\theta_{2}\left(\frac{v}{\tau}\mid-\frac{1}{\tau}\right) & =\sqrt{\frac{\tau}{i}}e^{\pi iv^{2}/\tau}\theta_{4}\left(v\mid\tau\right),\label{eq:1.30}\\
\theta_{3}\left(\frac{v}{\tau}\mid-\frac{1}{\tau}\right) & =\sqrt{\frac{\tau}{i}}e^{\pi iv^{2}/\tau}\theta_{3}\left(v\mid\tau\right),\label{eq:1.31}\\
\theta_{4}\left(\frac{v}{\tau}\mid-\frac{1}{\tau}\right) & =\sqrt{\frac{\tau}{i}}e^{\pi iv^{2}/\tau}\theta_{2}\left(v\mid\tau\right).\label{eq:1.32}\end{align}

\begin{lem}
\label{lem:2}For \begin{equation}
0<a<1,\quad n\in\mathbb{N},\quad\gamma>0,\label{eq:1.33}\end{equation}
 and\begin{equation}
q=e^{-2\pi\gamma^{-1}n^{-a}},\label{eq:1.34}\end{equation}
we have\begin{equation}
(q;q)_{\infty}=\sqrt{\gamma n^{a}}\exp\left\{ \frac{\pi}{12}\left((\gamma n^{a})^{-1}-\gamma n^{a}\right)\right\} \left\{ 1+\mathcal{O}\left(e^{-2\pi\gamma n^{a}}\right)\right\} ,\label{eq:1.35}\end{equation}
 and\begin{equation}
\frac{1}{(q;q)_{\infty}}=\frac{\exp\left\{ \frac{\pi}{12}\left(\gamma n^{a}-(\gamma n^{a})^{-1}\right)\right\} }{\sqrt{\gamma n^{a}}}\left\{ 1+\mathcal{O}\left(e^{-2\pi\gamma n^{a}}\right)\right\} \label{eq:1.36}\end{equation}
 as $n\to\infty$. 
\end{lem}
\begin{proof}
From formulas \eqref{eq:1.17}, \eqref{eq:1.19} and \eqref{eq:1.19}
we get \begin{align}
 & (q;q)_{\infty}=\exp\left(\pi\gamma^{-1}n^{-a}/12\right)\eta\left(\gamma^{-1}n^{-a}i\right)\label{eq:1.37}\\
 & =\sqrt{\gamma n^{a}}\exp\left(\pi\gamma^{-1}n^{-a}/12\right)\eta(\gamma n^{a}i)\nonumber \\
 & =\sqrt{\gamma n^{a}}\exp\left(\pi\gamma^{-1}n^{-a}/12-\pi\gamma n^{a}/12\right)\prod_{k=1}^{\infty}(1-e^{-2\pi\gamma kn^{a}})\nonumber \\
 & =\sqrt{\gamma n^{a}}\exp\left(\pi\gamma^{-1}n^{-a}/12-\pi\gamma n^{a}/12\right)\left\{ 1+\mathcal{O}\left(e^{-2\pi\gamma n^{a}}\right)\right\} ,\nonumber \end{align}
 and\begin{equation}
\frac{1}{(q;q)_{\infty}}=\frac{\exp\left(\pi\gamma n^{a}/12-\pi\gamma^{-1}n^{-a}/12\right)}{\sqrt{\gamma n^{a}}}\left\{ 1+\mathcal{O}\left(e^{-2\pi\gamma n^{a}}\right)\right\} \label{eq:1.38}\end{equation}
 as $n\to\infty$. 
\end{proof}
The Euler's $q$-Exponential is defined by \cite{Andrews4,Gasper,Ismail2,Koekoek}\begin{equation}
E_{q}(z):=(-z;q)_{\infty}=\sum_{k=0}^{\infty}\frac{q^{k(k-1)/2}}{(q;q)_{k}}z^{k},\quad z\in\mathbb{C}.\label{eq:1.39}\end{equation}
 The $q$-Gamma function is defined as \cite{Andrews4,Gasper,Ismail2,Koekoek}\begin{equation}
\Gamma_{q}(x):=\frac{(q;q)_{\infty}}{(q^{x};q)_{\infty}}(1-q)^{1-x}\quad x\in\mathbb{C}.\label{eq:1.40}\end{equation}
 Ramanujan's function $A_{q}(z)$ is defined as \cite{Ismail2}\begin{equation}
A_{q}(z):=\sum_{k=0}^{\infty}\frac{q^{k^{2}}}{(q;q)_{k}}(-z)^{k},\quad z\in\mathbb{C}.\label{eq:1.41}\end{equation}
Jackson's $q$-Bessel function of second kind \cite{Andrews4,Gasper,Ismail2,Koekoek}\begin{equation}
J_{\nu}^{(2)}(z;q):=\frac{(q^{\nu+1};q)_{\infty}}{(q;q)_{\infty}}\sum_{k=0}^{\infty}\frac{(-1)^{k}(z/2)^{\nu+2k}}{(q,q^{\nu+1};q)_{k}}q^{k(\nu+k)},\quad\nu>-1,\quad z\in\mathbb{C}.\label{eq:1.42}\end{equation}
Stieltjes\emph{-}Wigert orthogonal polynomials $\left\{ S_{n}(x;q)\right\} _{n=0}^{\infty}$
are defined as \cite{Ismail2}\begin{equation}
S_{n}(x;q):=\sum_{k=0}^{n}\frac{q^{k^{2}}(-x)^{k}}{(q;q)_{k}(q;q)_{n-k}},\quad x\in\mathbb{C},\quad n\in\mathbb{N}\cup\left\{ 0\right\} .\label{eq:1.43}\end{equation}
Stieltjes-Wigert orthogonal polynomials come from an indeterminant
moment problem. They satisfy the orthogonality relation \begin{equation}
\int_{0}^{\infty}S_{m}(x;q)S_{n}(x;q)w_{sw}(x)dx=\frac{q^{-n}}{(q;q)_{n}}\delta_{m,n},\quad n,m\in\mathbb{N}\cup\left\{ 0\right\} ,\label{eq:1.44}\end{equation}
 where\begin{align}
w_{sw}(x) & :=\sqrt{\frac{-1}{2\pi\log q}}\exp\left(\frac{1}{2\log q}\left[\log\left(\frac{x}{\sqrt{q}}\right)\right]^{2}\right),\quad x\in\mathbb{R}.\label{eq:1.45}\end{align}
 Clearly, the associated orthonormal Stieltjes-Wigert functions are
given by \begin{equation}
s_{n}(x;q):=\sqrt{q^{n}(q;q)_{n}w_{sw}(x)}S_{n}(x;q),\quad x\in\mathbb{R}^{+},\quad n\in\mathbb{N}\cup\left\{ 0\right\} .\label{eq:1.46}\end{equation}
The $q$-Laguerre orthogonal polynomials $\left\{ L_{n}^{(\alpha)}(x;q)\right\} _{n=0}^{\infty}$
are defined as \cite{Andrews4,Gasper,Ismail2,Koekoek} \begin{equation}
L_{n}^{(\alpha)}(x;q):=\sum_{k=0}^{n}\frac{q^{k^{2}+\alpha k}(-x)^{k}(q^{\alpha+1};q)_{n}}{(q;q)_{k}(q,q^{\alpha+1};q)_{n-k}},\quad\alpha>-1,\quad x\in\mathbb{C},\quad n\in\mathbb{N}\cup\left\{ 0\right\} .\label{eq:1.47}\end{equation}
The $q$-Laguerre orthogonal polynomials come from an indeterminate
moment problem. They satisfy the orthogonality relation \begin{equation}
\int_{0}^{\infty}L_{m}^{(\alpha)}(x;q)L_{n}^{(\alpha)}(x;q)w_{q\ell}(x)dx=\frac{(q^{\alpha+1};q)_{n}}{q^{n}(q;q)_{n}}\delta_{m,n},\quad\alpha>-1,\quad n,m\in\mathbb{N}\cup\left\{ 0\right\} ,\label{eq:1.48}\end{equation}
 where\begin{align}
w_{q\ell}(x) & :=-\frac{\sin(\pi\alpha)}{\pi}\frac{(q;q)_{\infty}}{(q^{-\alpha};q)_{\infty}}\frac{x^{\alpha}}{(-x;q)_{\infty}},\quad x\in\mathbb{R}^{+},\quad\alpha>-1.\label{eq:1.49}\end{align}
Clearly, the associated orthonormal $q$-Laguerre orthogonal functions
are given by \begin{equation}
\ell_{n}(x;q):=\sqrt{\frac{q^{n}(q;q)_{n}}{(q^{\alpha+1};q)_{n}}w_{q\ell}(x)}L_{n}^{(\alpha)}(x;q),\quad\alpha>-1,\quad n\in\mathbb{N}\cup\left\{ 0\right\} ,\quad x\in\mathbb{R}^{+}.\label{eq:1.50}\end{equation}
For any positive integer $n$, we define\begin{equation}
\chi(n):=2\left\{ \frac{n}{2}\right\} ,\label{eq:1.51}\end{equation}
then,\begin{equation}
\chi(n)=n-2\left\lfloor \frac{n}{2}\right\rfloor =\left\lfloor \frac{n+1}{2}\right\rfloor -\left\lfloor \frac{n}{2}\right\rfloor ,\label{eq:1.52}\end{equation}
where $\left\lfloor x\right\rfloor $ is the greatest integer less
than or equals to $x\in\mathbb{R}$ and $\left\{ x\right\} $ is the
fractional part of $x\in\mathbb{R}$.

\section{Main Results}

\subsection{Euler $q$-Exponential Function $E_{q}(z)$, $q$-Gamma Function
$\Gamma_{q}(x)$}

\begin{thm}
\label{thm:exponential}For \begin{equation}
0<a<\frac{1}{2},\quad n\in\mathbb{N},\quad u\in\mathbb{R},\quad q=\exp(-2n^{-a}\pi),\label{eq:2.1}\end{equation}
we have\begin{align}
E_{q}(\exp2\pi(u+n^{1-a}-\frac{1}{2}n^{-a}))= & \exp\left\{ \pi n^{-a}(n^{a}u+n)^{2}+\frac{\pi}{12}(n^{a}-n^{-a})\right\} \left\{ 1+\mathcal{O}\left(e^{-\pi n^{a}}\right)\right\} ,\label{eq:2.2}\end{align}
 and\begin{equation}
E_{q}(-\exp2\pi(u+n^{1-a}-\frac{1}{2}n^{-a}))=\frac{2\exp\left(\pi n^{-a}(n^{a}u+n)^{2}\right)\cos(\pi n^{a}u)}{(-1)^{n}\exp\frac{\pi}{12}(2n^{a}-n^{-a})}\left\{ 1+\mathcal{O}\left(e^{-2\pi n^{a}}\right)\right\} \label{eq:2.3}\end{equation}
 as $n\to\infty$. 
\end{thm}
Similarly, $\Gamma_{q}(z)$ has asymptotic behavior: 

\begin{thm}
\label{thm:gamma}For \begin{equation}
0<a<\frac{1}{2},\quad n\in\mathbb{N},\quad u\in\mathbb{R},\quad q=\exp(-2n^{-a}\pi),\label{eq:2.4}\end{equation}
we have\begin{align}
\frac{1}{\Gamma_{q}\left(\frac{1}{2}-n-n^{a}u\right)} & =\frac{2(-1)^{n}\exp\left(\pi n^{-a}(n^{a}u+n)^{2}\right)\cos(\pi n^{a}u)\left\{ 1+\mathcal{O}\left(e^{-2\pi n^{a}}\right)\right\} }{\sqrt{n^{a}}\exp\left(\pi n^{a}/12+\pi n^{-a}/6\right)\left(1-\exp(-2\pi n^{-a})\right)^{n+n^{a}u+1/2}},\label{eq:2.5}\end{align}
 and\begin{align}
\frac{1}{\Gamma_{q}\left(\frac{1}{2}+n+n^{a}u\right)} & =\frac{\exp(\pi n^{a}/12-\pi n^{-a}/12)\left\{ 1+\mathcal{O}\left(e^{-2\pi n^{a}}\right)\right\} }{\sqrt{n^{a}}(1-e^{-2\pi n^{-a}})^{1/2-n-n^{a}u}}\label{eq:2.6}\end{align}
 as $n\to\infty$. 
\end{thm}

\subsection{Ramanujan's Function $A_{q}(z)$}

\begin{thm}
\label{thm:ramanujan}For \begin{equation}
0<q<1,\quad z\in\mathbb{C}\backslash\left\{ 0\right\} ,\label{eq:2.7}\end{equation}
we have \begin{equation}
A_{q}(q^{-2n}z)=\frac{(-z)^{n}\left\{ \theta_{4}\left(z^{-1};q\right)+e(n)\right\} }{(q;q)_{\infty}q^{n^{2}}},\label{eq:2.8}\end{equation}
and\begin{equation}
|e(n)|\le4\theta_{3}\left(|z|^{-1};q\right)\left\{ \frac{q^{n/2}}{1-q}+\frac{q^{\left\lfloor n/2\right\rfloor ^{2}}}{\left|z\right|^{\left\lfloor n/2\right\rfloor }}\right\} \label{eq:2.9}\end{equation}
 for $n$ sufficiently large. 

Let\begin{equation}
q=\exp(-\pi n^{-a}),\quad0<a<\frac{1}{2},\quad n\in\mathbb{N},\quad u\in\mathbb{R},\label{eq:2.10}\end{equation}
then,\begin{equation}
A_{q}(-\exp2\pi(u+n^{1-a}))=\frac{\exp\left\{ \pi n^{-a}(n^{a}u+n)^{2}\right\} \left\{ 1+\mathcal{O}\left(\exp(-\pi n^{a})\right)\right\} }{\sqrt{2}\exp\left\{ \pi n^{-a}/24-\pi n^{a}/6\right\} },\label{eq:2.11}\end{equation}
and \begin{align}
A_{q}(\exp2\pi(u+n^{1-a})) & =\frac{\sqrt{2}\exp\left\{ \pi n^{-a}(n^{a}u+n)^{2}\right\} \left\{ \cos(\pi n^{a}u)+\mathcal{O}\left(\exp(-2\pi n^{a})\right)\right\} }{(-1)^{n}\exp\pi\left\{ n^{a}/12+n^{-a}/24\right\} }\label{eq:2.12}\end{align}
 as $n\to\infty$. 
\end{thm}

\subsection{Jackson's $q$-Bessel function of second kind $J_{\nu}^{(2)}(z;q)$}

\begin{thm}
\label{thm:jackson}For \begin{equation}
0<q<1,\quad z\in\mathbb{C}\backslash\left\{ 0\right\} ,\quad\nu>-1,\label{eq:2.13}\end{equation}
we have\begin{equation}
J_{\nu}^{(2)}\left(2\sqrt{zq^{-2n-\nu}};q\right)=\frac{z^{n+\nu/2}\left\{ \theta_{4}(z^{-1};q)+e(n)\right\} }{(-1)^{n}(q;q)_{\infty}^{2}q^{n^{2}+n\nu+\nu^{2}/2}},\label{eq:2.14}\end{equation}
 and\begin{equation}
|e(n)|\le12\theta_{3}(|z|^{-1};q)\left\{ \frac{q^{n/2}}{1-q}+\frac{q^{\left\lfloor n/2\right\rfloor ^{2}}}{|z|^{\left\lfloor n/2\right\rfloor }}\right\} \label{eq:2.15}\end{equation}
 for $n$ sufficiently large.

Let \begin{equation}
\nu>-1,\quad q=\exp(-\pi n^{-a}),\quad0<a<\frac{1}{2},\quad n\in\mathbb{N},\quad u\in\mathbb{R},\label{eq:2.16}\end{equation}
then, \begin{align}
 & J_{\nu}^{(2)}\left(2\imath\exp\left(\pi(u+n^{1-a}+\nu n^{-a}/2)\right);\exp(-\pi n^{-a})\right)\label{eq:2.17}\\
 & =\frac{\imath^{\nu}\exp\left\{ \pi n^{-a}\left(n^{a}u+n+\nu/2\right)^{2}\right\} \left\{ 1+\mathcal{O}\left(\exp(-\pi n^{a})\right)\right\} }{2\sqrt{n^{a}}\exp\pi\left\{ n^{-a}/12-n^{a}/3-n^{-a}\nu^{2}/4\right\} },\nonumber \end{align}
 and\begin{align}
 & J_{\nu}^{(2)}\left(2\exp\left(\pi(u+n^{1-a}+\nu n^{-a}/2)\right);\exp(-\pi n^{-a})\right)\label{eq:2.18}\\
 & =\frac{\exp\left\{ \pi n^{-a}\left(n^{a}u+n+\nu/2\right)^{2}\right\} \left\{ \cos(\pi n^{a}u)+\mathcal{O}\left(\exp(-2\pi n^{a})\right)\right\} }{(-1)^{n}\sqrt{n^{a}}\exp\pi\left\{ n^{-a}/12-n^{a}/12-\nu^{2}n^{-a}/4\right\} }\nonumber \end{align}
 as $n\to\infty$. 
\end{thm}

\subsection{Stieltjes-Wigert Orthogonal Polynomials $S_{n}(x;q)$ }

\begin{thm}
\label{thm:stieltjes-wigert}For\begin{equation}
0<q<1,\quad z\in\mathbb{C}\backslash\left\{ 0\right\} ,\label{eq:2.19}\end{equation}
 we have

\begin{equation}
S_{n}(zq^{-n};q)=\frac{(-z)^{\left\lfloor n/2\right\rfloor }\left\{ \theta_{4}(z^{-1}q^{\chi(n)};q)+e(n)\right\} }{(q;q)_{\infty}^{2}q^{\left\lfloor n/2\right\rfloor \left\lfloor (n+1)/2\right\rfloor }},\label{eq:2.20}\end{equation}
and\begin{align}
|e(n)| & \le12\theta_{3}\left(\left|z\right|^{-1}q^{\chi(n)};q\right)\left\{ \frac{q^{n/4}}{1-q}+|z|^{\left\lfloor n/4\right\rfloor }q^{\left\lfloor n/4\right\rfloor ^{2}-\chi(n)\left\lfloor n/4\right\rfloor }+\frac{q^{\left\lfloor n/4\right\rfloor ^{2}+\chi(n)\left\lfloor n/4\right\rfloor }}{|z|^{\left\lfloor n/4\right\rfloor }}\right\} \label{eq:2.21}\end{align}
 for $n$ sufficiently large. 

Let\begin{equation}
q=\exp(-2\pi n^{-a}),\quad0<a<\frac{1}{2},\quad n\in\mathbb{N},\quad u\in\mathbb{R},\label{eq:2.22}\end{equation}

then,\begin{align}
 & S_{n}(-\exp2\pi n^{-a}(n^{a}u+n);\exp(-2\pi n^{-a}))\label{eq:2.23}\\
 & =\frac{\exp\left\{ \pi n^{-a}(n^{a}u+n)^{2}/2\right\} \left\{ 1+\mathcal{O}\left(e^{-\pi n^{a}/2}\right)\right\} }{\sqrt{2n^{a}}\exp\left\{ \pi n^{-a}/6-\pi n^{a}/6\right\} },\label{eq:2.24}\end{align}
 and\begin{align}
 & S_{n}\left(\exp\left(2\pi n^{-a}(n^{a}u+n)\right);\exp(-2\pi n^{-a})\right)\label{eq:2.25}\\
 & =\sqrt{\frac{2}{n^{a}}}\frac{\exp\left\{ \frac{\pi n^{-a}}{2}(n^{a}u+n)^{2}\right\} \left\{ \cos\frac{\pi}{2}(n^{\alpha}u+n)+\mathcal{O}(e^{-\pi n^{a}})\right\} }{\exp\left\{ \pi n^{-a}/6-\pi n^{a}/24\right\} }\nonumber \end{align}
 for $n$ sufficiently large. Consequently,\begin{align}
 & s_{n}(\exp2\pi n^{-a}(n^{a}u+n);\exp(-2\pi n^{-a}))\label{eq:2.26}\\
 & =\frac{\exp(-u\pi/2)\left\{ \cos\frac{\pi}{2}\left(n^{\alpha}u+n\right)+\mathcal{O}\left(e^{-\pi n^{a}}\right)\right\} }{\sqrt{\pi}\exp(3\pi n^{1-a}/2+\pi n^{-a}/4)}\nonumber \end{align}

as $n\to\infty$. 
\end{thm}

\subsection{$q$-Laguerre Orthogonal Polynomials $L_{n}^{(\alpha)}(x;q)$}

\begin{thm}
\label{thm:laguerre}For\begin{equation}
0<q<1,\quad z\in\mathbb{C}\backslash\left\{ 0\right\} ,\quad\alpha>-1,\label{eq:2.27}\end{equation}
we have\begin{equation}
L_{n}^{(\alpha)}(zq^{-n-\alpha};q)=\frac{(-z)^{\left\lfloor n/2\right\rfloor }\left\{ \theta_{4}(z^{-1}q^{\chi(n)};q)+e(n)\right\} }{(q;q)_{\infty}^{2}q^{\left\lfloor n/2\right\rfloor \left\lfloor (n+1)/2\right\rfloor }},\label{eq:2.28}\end{equation}
and\begin{equation}
|e(n)|\le60\theta_{3}(|z|^{-1}q^{\chi(n)};q)\left\{ \frac{q^{\left\lfloor n/4\right\rfloor ^{2}+\chi(n)\left\lfloor n/4\right\rfloor }}{|z|^{\left\lfloor n/4\right\rfloor }}+|z|^{\left\lfloor n/4\right\rfloor }q^{\left\lfloor n/4\right\rfloor ^{2}-\chi(n)\left\lfloor n/4\right\rfloor }+\frac{q^{n/4}}{1-q}\right\} \label{eq:2.29}\end{equation}
for sufficiently $n$. Let\begin{equation}
q=\exp(-2\pi n^{-a}),\quad0<a<\frac{1}{2},\quad\alpha>-1,\quad n\in\mathbb{N},\quad u\in\mathbb{R},\label{eq:2.30}\end{equation}
then, \begin{align}
 & L_{n}^{\alpha}(-\exp2\pi(u+n^{1-a}+\alpha n^{-a});\exp(-2\pi n^{-a}))\label{eq:2.31}\\
 & =\frac{\exp\left\{ \frac{\pi n^{-a}}{2}(n^{a}u+n)^{2}\right\} \left\{ 1+\mathcal{O}\left(\exp(-\pi n^{a}/2)\right)\right\} }{\sqrt{2n^{a}}\exp\pi\left\{ \pi n^{-a}/6-\pi n^{a}/6\right\} },\nonumber \end{align}
 and\begin{align}
 & L_{n}^{\alpha}(\exp2\pi(u+n^{1-a}+\alpha n^{-a});\exp(-2\pi n^{-a}))\label{eq:2.32}\\
 & =\sqrt{\frac{2}{n^{a}}}\frac{\exp\left\{ \frac{\pi n^{-a}}{2}\left(n^{a}u+n\right)^{2}\right\} \left\{ \cos\frac{\pi}{2}(n^{a}u+n)+\mathcal{O}\left(\exp(-\pi n^{a})\right)\right\} }{\exp\left\{ \pi n^{a}/24+\pi n^{-a}/6\right\} }\nonumber \end{align}

as $n\to\infty$. Consequently,\begin{align}
 & \ell_{n}^{(\alpha)}(\exp2\pi(u+n^{1-a}+\alpha n^{-a});\exp(-2\pi n^{-a}))\label{eq:2.33}\\
 & =\frac{1}{\sqrt{\pi}}\frac{\exp(-\pi u/2)\left\{ \cos\frac{\pi}{2}(n^{a}u+n)+\mathcal{O}\left(\exp(-\pi n^{a})\right)\right\} }{\exp\frac{\pi}{2}(3n^{1-a}+\alpha^{2}n^{a}+n^{a}/6+2\alpha n^{-a}+n^{-a}/2-\alpha^{2}n^{-a})}\nonumber \end{align}
 as $n\to\infty$.
\end{thm}

\section{Proofs }

\subsection{Proof for Theorem \ref{thm:exponential} }

\begin{proof}
It is clear that\begin{align}
E_{q}(zq^{1/2-n}) & =(-zq^{1/2-n};q)_{\infty}=\frac{z^{n}q^{-n^{2}/2}(q,-zq^{1/2},-z^{-1}q^{1/2};q)_{\infty}}{(q,-z^{-1}q^{1/2+n};q)_{\infty}}.\label{eq:3.1}\end{align}
Thus,\begin{align}
 & E_{q}(q^{-n+1/2}e^{2\pi u})=(-q^{-n+1/2}e^{2\pi u};q)_{\infty}\label{eq:3.2}\\
 & =\frac{q^{-n^{2}/2}e^{2\pi nu}\left(q,-q^{1/2}e^{2\pi u},-q^{1/2}e^{-2\pi u};q\right)_{\infty}}{\left(q,-q^{n+1/2}e^{-2\pi u};q\right)_{\infty}}\nonumber \\
 & =\frac{q^{-n^{2}/2}e^{2\pi nu}\theta_{3}(e^{2\pi u};q^{1/2})}{\left(q,-q^{n+1/2}e^{-2\pi u};q\right)_{\infty}},\nonumber \end{align}
and\begin{align}
 & E_{q}(-q^{-n+1/2}e^{2\pi u})=(q^{-n+1/2}e^{2\pi u};q)_{\infty}\label{eq:3.3}\\
 & =\frac{(-1)^{n}q^{-n^{2}/2}e^{2\pi nu}\left(q,q^{1/2}e^{2\pi u},q^{1/2}e^{-2\pi u};q\right)_{\infty}}{\left(q,q^{n+1/2}e^{-2\pi u};q\right)_{\infty}}\nonumber \\
 & =\frac{(-1)^{n}q^{-n^{2}/2}e^{2\pi nu}\theta_{4}(e^{2\pi u};q^{1/2})}{\left(q,q^{n+1/2}e^{-2\pi u};q\right)_{\infty}}.\nonumber \end{align}
Since \begin{equation}
\frac{1}{(q;q)_{\infty}}=\frac{\exp\left\{ \frac{\pi}{12}\left(n^{a}-n^{-a}\right)\right\} }{\sqrt{n^{a}}}\left\{ 1+\mathcal{O}\left(e^{-2\pi n^{a}}\right)\right\} ,\label{eq:3.4}\end{equation}
\begin{equation}
\frac{1}{\left(-q^{n+1/2}e^{-2\pi u};q\right)_{\infty}}=1+\mathcal{O}\left(n^{a}e^{-2\pi n^{1-a}}\right),\label{eq:3.5}\end{equation}
\begin{equation}
\frac{1}{\left(q^{n+1/2}e^{-2\pi u};q\right)_{\infty}}=1+\mathcal{O}\left(n^{a}e^{-2\pi n^{1-a}}\right),\label{eq:3.6}\end{equation}
\begin{align}
 & \theta_{4}(e^{2\pi u};q^{1/2})=\theta_{4}(ui\vert n^{-a}i)=\sqrt{n^{a}}\exp(\pi n^{a}u^{2})\theta_{2}(n^{a}u\vert n^{a}i)\label{eq:3.7}\\
 & =2\sqrt{n^{a}}\exp\left(\pi n^{a}u^{2}-\frac{\pi n^{a}}{4}\right)\cos(\pi n^{a}u)\left\{ 1+\mathcal{O}\left(e^{-2\pi n^{a}}\right)\right\} ,\nonumber \end{align}
and\begin{align}
 & \theta_{3}(e^{2\pi u};q^{1/2})=\theta_{3}(ui\vert n^{-a}i)=\sqrt{n^{a}}\exp(\pi n^{a}u^{2})\theta_{3}(n^{a}u\vert n^{a}i)\label{eq:3.8}\\
 & =\sqrt{n^{a}}\exp(\pi n^{a}u^{2})\left\{ 1+\mathcal{O}\left(e^{-\pi n^{a}}\right)\right\} \nonumber \end{align}
as $n\to\infty$. Thus, \begin{equation}
E_{q}(\exp2\pi(u+n^{1-a}-\frac{1}{2}n^{-a}))=\exp\left\{ \pi n^{-a}(n^{a}u+n)^{2}+\frac{\pi}{12}(n^{a}-n^{-a})\right\} \left\{ 1+\mathcal{O}\left(e^{-\pi n^{a}}\right)\right\} ,\label{eq:3.9}\end{equation}
 and

\begin{equation}
E_{q}(-\exp2\pi(u+n^{1-a}-\frac{1}{2}n^{-a}))=\frac{2\exp(\pi n^{-a}(n^{a}u+n)^{2})\cos(\pi n^{a}u)}{(-1)^{n}\exp\frac{\pi}{12}(2n^{a}-n^{-a})}\left\{ 1+\mathcal{O}\left(e^{-2\pi n^{a}}\right)\right\} \label{eq:3.10}\end{equation}
 as $n\to\infty$. 
\end{proof}

\subsection{Proof for Theorem \ref{thm:gamma}}

\begin{proof}
Observe that \begin{align}
 & \frac{1}{\Gamma_{q}\left(1/2-n-n^{a}u\right)}=\frac{(q^{1/2-n}e^{2\pi u};q)_{\infty}}{(q;q)_{\infty}(1-q)^{n+n^{a}u+1/2}}\label{eq:3.11}\\
 & =\frac{(q,q^{1/2}e^{-2\pi u},q^{1/2}e^{2\pi u};q)_{\infty}q^{-n^{2}/2}e^{2\pi nu}}{(-1)^{n}(1-q)^{n+n^{a}u+1/2}(q,q,q^{n+1/2}e^{-2\pi u};q)_{\infty}}.\nonumber \end{align}
 Since\begin{align}
\frac{1}{(q,q,q^{n+1/2}e^{-2\pi u};q)_{\infty}} & =n^{-a}\exp\left(\pi n^{a}/6-\pi n^{-a}/6\right)\left\{ 1+\mathcal{O}\left(e^{-2\pi n^{a}}\right)\right\} ,\label{eq:3.12}\end{align}
and \begin{align}
 & (q,q^{1/2}e^{-2\pi u},q^{1/2}e^{2\pi u};q)_{\infty}=\theta_{4}(ui\mid n^{-a}i)=n^{a/2}e^{\pi n^{a}u^{2}}\theta_{2}(n^{a}u\mid n^{a}i)\label{eq:3.13}\\
 & =2n^{a/2}\exp\pi n^{a}(u^{2}-1/4)\cos(n^{a}u\pi)\left\{ 1+\mathcal{O}\left(e^{-2\pi n^{a}}\right)\right\} \nonumber \end{align}
 as $n\to\infty$. Thus,\begin{align}
\frac{1}{\Gamma_{q}\left(1/2-n-n^{a}u\right)} & =\frac{2\exp\pi\left(n^{-a}(n^{a}u+n)^{2}\right)\cos(\pi n^{a}u)\left\{ 1+\mathcal{O}\left(e^{-2\pi n^{a}}\right)\right\} }{(-1)^{n}n^{a/2}\left(1-e^{-2\pi n^{-a}}\right)^{n+n^{a}u+1/2}\exp\left(\pi n^{a}/12+\pi n^{-a}/6\right)}\label{eq:3.14}\end{align}
 as $n\to\infty$.

Similarly we could prove that 

\begin{align}
\frac{1}{\Gamma_{q}(1/2+n+n^{a}u)} & =\frac{\exp(\pi n^{a}/12-\pi n^{-a}/12)\left\{ 1+\mathcal{O}\left(e^{-2\pi n^{a}}\right)\right\} }{n^{a/2}(1-e^{-2\pi n^{-a}})^{1/2-n-n^{a}u}}\label{eq:3.15}\end{align}
 as $n\to\infty$. 
\end{proof}

\subsection{Proof for Theorem \ref{thm:ramanujan}}

\begin{proof}
Write\begin{align}
A_{q}(q^{-2n}z)(q;q)_{\infty} & =\sum_{k=0}^{\infty}\frac{(q;q)_{\infty}}{(q;q)_{k}}q^{k^{2}}\left(-zq^{-2n}\right)^{k}\label{eq:3.16}\\
 & =\sum_{k=0}^{n}f(k)q^{k^{2}}\left(-zq^{-2n}\right)^{k}+\sum_{k=n+1}^{\infty}f(k)q^{k^{2}}\left(-zq^{-2n}\right)^{k}\nonumber \\
 & =s_{1}+s_{2},\nonumber \end{align}
where\begin{equation}
f(k)=(q^{k+1};q)_{\infty},\quad k\in\mathbb{N}\cup\left\{ 0\right\} .\label{eq:3.17}\end{equation}
Clearly, \begin{equation}
|f(k)|\le1,\quad k\in\mathbb{N}\cup\left\{ 0\right\} .\label{eq:3.18}\end{equation}
 Reverse summation order in $s_{1}$ to get\begin{align}
\frac{s_{1}q^{n^{2}}}{(-z)^{n}} & =\sum_{k=0}^{n}f(n-k)q^{k^{2}}(-z^{-1})^{k}.\label{eq:3.19}\end{align}
 By Lemma\ref{lem:1} we have

\begin{equation}
|f(n-k)-1|=\left|r_{1}\left(q;n-k\right)\right|\le\frac{2q^{n/2}}{1-q}\label{eq:3.20}\end{equation}
 for $0\le k\le\left\lfloor \frac{n}{2}\right\rfloor $ and $n$ sufficiently
large. Hence,

\begin{align}
\frac{s_{1}q^{n^{2}}}{(-z)^{n}} & =\sum_{k=0}^{n}f(n-k)q^{k^{2}}(-z^{-1})^{k}\label{eq:3.21}\\
 & =\sum_{k=0}^{\infty}q^{k^{2}}(-z^{-1})^{k}-\sum_{k=\left\lfloor n/2\right\rfloor }^{\infty}q^{k^{2}}(-z^{-1})^{k}\nonumber \\
 & +\sum_{k=0}^{\left\lfloor n/2\right\rfloor -1}q^{k^{2}}(-z^{-1})^{k}\left\{ f(n-k)-1\right\} +\sum_{k=\left\lfloor n/2\right\rfloor }^{n}q^{k^{2}}(-z^{-1})^{k}f(n-k)\nonumber \\
 & =\sum_{k=0}^{\infty}q^{k^{2}}(-q^{\chi(n)}z^{-1})^{k}+s_{11}+s_{12}+s_{13}.\nonumber \end{align}
 Thus,

\begin{align}
|s_{11}+s_{13}| & \le2\sum_{k=\left\lfloor n/2\right\rfloor }^{\infty}q^{k^{2}}\left(\frac{1}{|z|}\right)^{k}\le2\frac{q^{\left\lfloor n/2\right\rfloor ^{2}}}{\left|z\right|^{\left\lfloor n/2\right\rfloor }}\sum_{k=0}^{\infty}q^{k^{2}}\left(\frac{1}{|z|}\right)^{k}\le2\theta_{3}\left(|z|^{-1};q\right)\frac{q^{\left\lfloor n/2\right\rfloor ^{2}}}{\left|z\right|^{\left\lfloor n/2\right\rfloor }}\label{eq:3.22}\end{align}
 and\begin{align}
|s_{12}| & \le\frac{2q^{n/2}}{1-q}\sum_{k=0}^{\infty}q^{k^{2}}\left(\frac{1}{|z|}\right)^{k}\le\frac{2\theta_{3}\left(|z|^{-1};q\right)}{1-q}q^{n/2}\label{eq:3.23}\end{align}
 for $n$ sufficiently large. Let \begin{equation}
e_{1}(n)=s_{11}+s_{12}+s_{13},\label{eq:3.24}\end{equation}
 then,\begin{equation}
|e_{1}(n)|\le2\theta_{3}\left(|z|^{-1};q\right)\left\{ \frac{q^{n/2}}{1-q}+\frac{q^{\left\lfloor n/2\right\rfloor ^{2}}}{\left|z\right|^{\left\lfloor n/2\right\rfloor }}\right\} \label{eq:3.25}\end{equation}
 and\begin{equation}
\frac{s_{1}q^{n^{2}}}{(-z)^{n}}=\sum_{k=0}^{\infty}q^{k^{2}}(-z^{-1})^{k}+e_{1}(n),\label{eq:3.26}\end{equation}
 for $n$ sufficiently large.

Shift the summation index from $k$ to $k+n$ in $s_{2}$ we obtain\begin{align}
\frac{s_{2}q^{n^{2}}}{(-z)^{n}} & =\sum_{k=1}^{\infty}q^{k^{2}}(-z)^{k}f(n+k).\label{eq:3.27}\end{align}
Lemma\ref{lem:1} implies that \begin{equation}
|f(n+k)-1|=\left|r_{1}\left(q;n+k\right)\right|\le\frac{2q^{n}}{1-q}\label{eq:3.28}\end{equation}
 for $k\in\mathbb{N}$ and $n$ sufficiently large. Thus,

\begin{align}
\frac{s_{2}q^{n^{2}}}{(-z)^{n}} & =\sum_{k=1}^{\infty}q^{k^{2}}(-z)^{k}f(n+k)\label{eq:3.29}\\
 & =\sum_{k=1}^{\infty}q^{k^{2}}(-z)^{k}+\sum_{k=1}^{\infty}q^{k^{2}}(-z)^{k}\left\{ f(n+k)-1\right\} \nonumber \\
 & =\sum_{k=-1}^{-\infty}q^{k^{2}}(-z^{-1})^{k}+e_{2}(n),\nonumber \end{align}
 and\begin{align}
|e_{2}(n)| & \le\frac{2q^{n}}{1-q}\sum_{k=1}^{\infty}q^{k^{2}}|z|^{k}\le\frac{2\theta_{3}\left(|z|^{-1};q\right)}{1-q}q^{n/2}\label{eq:3.30}\end{align}
 for $n$ sufficiently large. Hence,\begin{equation}
\frac{A_{q}(q^{-2n}z)(q;q)_{\infty}q^{n^{2}}}{(-z)^{n}}=\theta_{4}\left(z^{-1};q\right)+e(n),\label{eq:3.31}\end{equation}
 with\begin{equation}
|e(n)|\le4\theta_{3}\left(|z|^{-1};q\right)\left\{ \frac{q^{n/2}}{1-q}+\frac{q^{\left\lfloor n/2\right\rfloor ^{2}}}{\left|z\right|^{\left\lfloor n/2\right\rfloor }}\right\} \label{eq:3.32}\end{equation}
 for $n$ sufficiently large.

From the properties of $\theta_{3}$ and $\theta_{4}$ we get

\begin{align}
\theta_{3}\left(\left|z\right|^{-1};q\right) & =\theta_{3}\left(iu\mid n^{-a}i\right)\label{eq:3.33}\\
 & =\sqrt{n^{a}}\exp\left\{ \pi n^{a}u^{2}\right\} \theta_{3}\left(n^{a}u\mid n^{a}i\right)\nonumber \\
 & =\sqrt{n^{a}}\exp\left\{ \pi n^{a}u^{2}\right\} \left\{ 1+\mathcal{O}\left(e^{-\pi n^{a}}\right)\right\} ,\nonumber \end{align}
and \begin{align}
\theta_{4}\left(z^{-1};q\right) & =\theta_{4}\left(iu\mid n^{-a}i\right)\label{eq:3.34}\\
 & =\sqrt{n^{a}}\exp\left\{ \pi n^{a}u^{2}\right\} \theta_{2}\left(n^{a}u\mid n^{a}i\right)\nonumber \\
 & =2\sqrt{n^{a}}\exp\left\{ \pi n^{a}u^{2}-\pi n^{a}/4\right\} \cos(\pi n^{a}u)\left\{ 1+\mathcal{O}\left(\exp(-2\pi n^{a})\right)\right\} \nonumber \end{align}
 as $n\to\infty$. Observe that\begin{equation}
n^{a}\ll n^{1-a}\label{eq:3.35}\end{equation}
 for $0<a<\frac{1}{2}$ as $n\to\infty$, thus, \begin{equation}
A_{q}(-\exp2\pi(u+n^{1-a}))=\frac{\exp\left\{ \pi n^{-a}(n^{a}u+n)^{2}\right\} \left\{ 1+\mathcal{O}\left(\exp(-\pi n^{a})\right)\right\} }{\sqrt{2}\exp\left\{ \pi n^{-a}/24-\pi n^{a}/6\right\} },\label{eq:3.36}\end{equation}
and \begin{align}
A_{q}(\exp2\pi(u+n^{1-a})) & =\frac{\sqrt{2}\exp\left\{ \pi n^{-a}(n^{a}u+n)^{2}\right\} \left\{ \cos(\pi n^{a}u)+\mathcal{O}\left(\exp(-2\pi n^{a})\right)\right\} }{(-1)^{n}\exp\pi\left\{ n^{a}/12+n^{-a}/24\right\} }\label{eq:3.37}\end{align}
 as $n\to\infty$. 
\end{proof}

\subsection{Proof for Theorem \ref{thm:jackson} }

\begin{proof}
It is clear that\begin{equation}
\frac{J_{\nu}^{(2)}\left(2\sqrt{xq^{-\nu}};q\right)(q;q)_{\infty}^{2}}{(xq^{-\nu})^{\nu/2}}=\sum_{k=0}^{\infty}q^{k^{2}}(-x)^{k}f(k),\label{eq:3.38}\end{equation}
where\begin{equation}
f(k):=(q^{k+1},q^{\nu+1+k};q)_{\infty},\quad k\in\mathbb{N}\cup\left\{ 0\right\} .\label{eq:3.9}\end{equation}
For $\nu>-1$, we have\begin{equation}
|f(k)|\le1,\quad k\in\mathbb{N}\cup\left\{ 0\right\} .\label{eq:3.40}\end{equation}
It is clear that\begin{align}
\frac{J_{\nu}^{(2)}\left(2\sqrt{zq^{-2n-\nu}};q\right)(q;q)_{\infty}^{2}q^{n^{2}+n\nu+\nu^{2}/2}}{(-1)^{n}z^{n+\nu/2}} & =\sum_{k=0}^{n}q^{k^{2}}(-z^{-1})^{k}f(n-k)+\sum_{k=1}^{\infty}q^{k^{2}}(-z)^{k}f(n+k).\label{eq:3.41}\end{align}
Notice that for $k\ge0$, \begin{align}
f(k)-1 & =\left\{ r_{1}\left(q;k\right)+1\right\} \left\{ r_{1}\left(q^{\nu+1};k\right)+1\right\} -1=r_{1}(q;k)+r_{1}(q^{\nu+1};k)+r_{1}(q;k)r_{1}(q^{\nu+1};k),\label{eq:3.42}\end{align}
 then, for sufficiently large $n$, \begin{equation}
0<\frac{2q^{n/2}}{1-q}<1,\label{eq:3.43}\end{equation}
and Lemma\ref{lem:1} implies

\begin{equation}
|f(n-k)-1|\le\frac{6q^{n/2}}{1-q},\quad0\le k\le\left\lfloor \frac{n}{2}\right\rfloor ,\label{eq:3.44}\end{equation}
and\begin{equation}
|f(n+k)-1|\le\frac{6q^{n}}{1-q}<\frac{6q^{n/2}}{1-q},\quad k\in\mathbb{N}.\label{eq:3.45}\end{equation}
Similar to the proof of Theorem \ref{thm:ramanujan}, we have\begin{equation}
J_{\nu}^{(2)}\left(2\sqrt{zq^{-2n-\nu}};q\right)=\frac{z^{n+\nu/2}\left\{ \theta_{4}(z^{-1};q)+e(n)\right\} }{(-1)^{n}(q;q)_{\infty}^{2}q^{n^{2}+n\nu+\nu^{2}/2}},\label{eq:3.46}\end{equation}
 and\begin{equation}
|e(n)|\le12\theta_{3}(|z|^{-1};q)\left\{ \frac{q^{n/2}}{1-q}+\frac{q^{\left\lfloor n/2\right\rfloor ^{2}}}{|z|^{\left\lfloor n/2\right\rfloor }}\right\} \label{eq:3.47}\end{equation}
 for $n$ sufficiently large. The rest of the proof is very similar
to the corresponding part for Theorem \ref{thm:ramanujan}.
\end{proof}

\subsection{Proof for Theorem \ref{thm:stieltjes-wigert}}

\begin{proof}
As in the proof for Theorem \ref{thm:ramanujan} we have\begin{align}
\frac{(q;q)_{\infty}^{2}S_{n}(zq^{-n};q)}{(-z)^{\left\lfloor n/2\right\rfloor }q^{-\left\lfloor n/2\right\rfloor \left\lfloor (n+1)/2\right\rfloor }} & =\sum_{k=0}^{\left\lfloor \frac{n}{2}\right\rfloor }q^{k^{2}}(-z^{-1}q^{\chi(n)})^{k}f(\left\lfloor \frac{n}{2}\right\rfloor -k)+\sum_{k=0}^{\left\lfloor \frac{n+1}{2}\right\rfloor }q^{k^{2}}(-zq^{-\chi(n)})^{k}f(\left\lfloor \frac{n}{2}\right\rfloor +k),\label{eq:3.48}\end{align}
where\begin{equation}
f(k)=(q^{k+1};q)_{\infty}(q^{n-k+1};q)_{\infty},\quad k=0,1...,n.\label{eq:3.49}\end{equation}
Clearly\begin{equation}
|f(k)|\le1,\quad k=0,1,...,n.\label{eq:3.50}\end{equation}
Observe that\begin{equation}
f(k)-1=r_{1}(q;k)+r_{1}(q;n-k)+r_{1}(q;k)r_{1}(q;n-k),\label{eq:3.51}\end{equation}
and if \begin{equation}
0<\frac{2q^{n/4}}{1-q}<1\label{eq:3.52}\end{equation}
for $n$ sufficiently large, then Lemma\ref{lem:1} gives \begin{equation}
|f(\left\lfloor \frac{n}{2}\right\rfloor -k)-1|\le\frac{6q^{n/4}}{1-q},\quad0\le k\le\left\lfloor \frac{n}{4}\right\rfloor -1,\label{eq:3.53}\end{equation}
and\begin{equation}
|f(\left\lfloor \frac{n}{2}\right\rfloor +k)-1|\le\frac{6q^{n/4}}{1-q},\quad k\in\mathbb{N}.\label{eq:3.54}\end{equation}
 Thus, 

\begin{equation}
\frac{(q;q)_{\infty}^{2}S_{n}(zq^{-n};q)}{(-z)^{\left\lfloor n/2\right\rfloor }q^{-\left\lfloor n/2\right\rfloor \left\lfloor (n+1)/2\right\rfloor }}=\theta_{4}(z^{-1}q^{\chi(n)};q)+e(n),\label{eq:3.55}\end{equation}
 where\begin{align}
e(n) & =-\sum_{k=\left\lfloor n/4\right\rfloor }^{\infty}q^{k^{2}}(-z^{-1}q^{\chi(n)})^{k}+\sum_{k=0}^{\left\lfloor n/4\right\rfloor -1}q^{k^{2}}(-z^{-1}q^{\chi(n)})^{k}\left(f(\left\lfloor \frac{n}{2}\right\rfloor -k)-1\right)\label{eq:3.56}\\
 & +\sum_{k=\left\lfloor n/4\right\rfloor }^{\left\lfloor n/2\right\rfloor }q^{k^{2}}(-z^{-1}q^{\chi(n)})^{k}f(\left\lfloor \frac{n}{2}\right\rfloor -k)-\sum_{k=\left\lfloor n/4\right\rfloor }^{\infty}q^{k^{2}}(-zq^{-\chi(n)})^{k}\nonumber \\
 & +\sum_{k=1}^{\left\lfloor n/4\right\rfloor -1}q^{k^{2}}(-zq^{-\chi(n)})^{k}\left(f(\left\lfloor \frac{n}{2}\right\rfloor +k)-1\right)+\sum_{k=\left\lfloor n/4\right\rfloor }^{\left\lfloor (n+1)/2\right\rfloor }q^{k^{2}}(-zq^{-\chi(n)})^{k}f(\left\lfloor \frac{n}{2}\right\rfloor +k).\nonumber \end{align}
Therefore,\begin{equation}
S_{n}(zq^{-n};q)=\frac{(-z)^{\left\lfloor n/2\right\rfloor }\left\{ \theta_{4}(z^{-1}q^{\chi(n)};q)+e(n)\right\} }{(q;q)_{\infty}^{2}q^{\left\lfloor n/2\right\rfloor \left\lfloor (n+1)/2\right\rfloor }},\label{eq:3.57}\end{equation}
and\begin{align}
|e(n)| & \le12\theta_{3}\left(\left|z\right|^{-1}q^{\chi(n)};q\right)\left\{ \frac{q^{n/4}}{1-q}+|z|^{\left\lfloor n/4\right\rfloor }q^{\left\lfloor n/4\right\rfloor ^{2}-\chi(n)\left\lfloor n/4\right\rfloor }+\frac{q^{\left\lfloor n/4\right\rfloor ^{2}+\chi(n)\left\lfloor n/4\right\rfloor }}{|z|^{\left\lfloor n/4\right\rfloor }}\right\} \label{eq:3.58}\end{align}
 for $n$ sufficiently large. 

From the relations \begin{equation}
\left\lfloor \frac{n}{2}\right\rfloor =\frac{n-\chi(n)}{2},\quad\left\lfloor \frac{n+1}{2}\right\rfloor =\frac{n+\chi(n)}{2},\quad n\in\mathbb{N},\label{eq:3.59}\end{equation}
Lemma \ref{lem:2} and transformation formulas of theta functions
we get\begin{align}
 & S_{n}(-\exp2\pi n^{-a}(n^{a}u+n);\exp(-2\pi n^{-a}))\label{eq:3.60}\\
 & =\frac{\exp\left\{ \pi n^{-a}(n^{a}u+n)^{2}/2\right\} \left\{ 1+\mathcal{O}\left(e^{-\pi n^{a}/2}\right)\right\} }{\sqrt{2n^{a}}\exp\left\{ \pi n^{-a}/6-\pi n^{a}/6\right\} },\nonumber \end{align}
and 

\begin{align}
 & S_{n}(\exp2\pi n^{-a}(n^{a}u+n);\exp(-2\pi n^{-a}))\label{eq:3.61}\\
 & =\sqrt{\frac{2}{n^{a}}}\frac{\exp\left\{ \frac{\pi n^{-a}}{2}(n^{a}u+n)\right\} \left\{ \cos\frac{\pi}{2}\left(n^{\alpha}u+n\right)+\mathcal{O}\left(e^{-\pi n^{a}}\right)\right\} }{\exp\left\{ \pi n^{-a}/6-\pi n^{a}/24\right\} }\nonumber \end{align}
 as $n\to\infty$.

Lemma\ref{lem:1} and Lemma\ref{lem:2} also give us that \begin{align}
(q;q)_{n} & =\frac{(q;q)_{\infty}}{(q^{n+1};q)_{\infty}}=\sqrt{n^{a}}\exp\left(\pi n^{-a}/12-\pi n^{a}/12\right)\left\{ 1+\mathcal{O}\left(e^{-2\pi n^{a}}\right)\right\} ,\label{eq:3.62}\end{align}
 hence\begin{equation}
\sqrt{w_{sw}(\exp2\pi(u+n^{1-a}))}=\frac{\sqrt[4]{n^{a}}}{\sqrt{2\pi}}\exp\left(-\frac{n^{a}\pi}{2}(u+n^{1-a}+n^{-a}/2)^{2}\right),\label{eq:3.63}\end{equation}
 thus\begin{align}
 & \sqrt{q^{n}(q;q)_{n}w_{sw}(\exp2\pi(u+n^{1-a}))}\label{eq:3.64}\\
 & =\sqrt{\frac{n^{a}}{2\pi}}\exp\left(-\frac{n^{-a}\pi}{2}(n^{a}u+n)^{2}-\frac{u\pi}{2}\right)\nonumber \\
 & \times\exp\left(-\frac{\pi n^{-a}}{12}-\frac{\pi n^{a}}{24}-\frac{3\pi n^{1-a}}{2}\right)\left\{ 1+\mathcal{O}\left(e^{-2\pi n^{a}}\right)\right\} \nonumber \end{align}
 as $n\to\infty$. Then we get \begin{align}
 & s_{n}(\exp2\pi n^{-a}(n^{a}u+n);\exp(-2\pi n^{-a}))\label{eq:3.65}\\
 & =\frac{\exp(-u\pi/2)\left\{ \cos\frac{\pi}{2}\left(n^{\alpha}u+n\right)+\mathcal{O}\left(e^{-\pi n^{a}}\right)\right\} }{\sqrt{\pi}\exp(3\pi n^{1-a}/2+\pi n^{-a}/4)}\nonumber \end{align}
 as $n\to\infty$. 
\end{proof}

\subsection{Proof for Theorem \ref{thm:laguerre}}

\begin{proof}
Write\begin{align}
\frac{L_{n}^{(\alpha)}(zq^{-n-\alpha};q)(q;q)_{\infty}^{2}}{(-z)^{\left\lfloor n/2\right\rfloor }q^{-\left\lfloor n/2\right\rfloor \left\lfloor (n+1)/2\right\rfloor }} & =\sum_{k=0}^{\left\lfloor n/2\right\rfloor }q^{k^{2}}(-z^{-1}q^{\chi(n)})^{k}g(\left\lfloor \frac{n}{2}\right\rfloor -k)+\sum_{k=0}^{\left\lfloor (n+1)/2\right\rfloor }q^{k^{2}}(-zq^{-\chi(n)})^{k}g(\left\lfloor \frac{n}{2}\right\rfloor +k),\label{eq:3.66}\end{align}
where \begin{equation}
g(k)=\frac{(q^{k+1};q)_{\infty}(q^{n-k+1};q)_{\infty}(q^{\alpha+n-k+1};q)_{\infty}}{(q^{\alpha+n+1};q)_{\infty}},\quad k=0,1,\dotsc,n.\label{eq:3.67}\end{equation}
Then\begin{equation}
|g(k)|\le1,\quad k=0,1,...,n.\label{eq:3.68}\end{equation}
Thus,\begin{equation}
\frac{L_{n}^{(\alpha)}(zq^{-n-\alpha};q)(q;q)_{\infty}^{2}}{(-z)^{\left\lfloor n/2\right\rfloor }q^{-\left\lfloor n/2\right\rfloor \left\lfloor (n+1)/2\right\rfloor }}=\theta_{4}(z^{-1}q^{\chi(n)};q)+e(n),\label{eq:3.69}\end{equation}
where \begin{align}
e(n) & =-\sum_{k=\left\lfloor n/4\right\rfloor }^{\infty}q^{k^{2}}(-z^{-1}q^{\chi(n)})^{k}+\sum_{k=0}^{\left\lfloor n/4\right\rfloor -1}q^{k^{2}}(-z^{-1}q^{\chi(n)})^{k}(g(\left\lfloor \frac{n}{2}\right\rfloor -k)-1)\label{eq:3.70}\\
 & +\sum_{k=\left\lfloor n/4\right\rfloor }^{\left\lfloor n/2\right\rfloor }q^{k^{2}}(-z^{-1}q^{\chi(n)})^{k}g(\left\lfloor \frac{n}{2}\right\rfloor -k)-\sum_{k=\left\lfloor n/4\right\rfloor }^{\infty}q^{k^{2}}(-zq^{-\chi(n)})^{k}\nonumber \\
 & +\sum_{k=0}^{\left\lfloor n/4\right\rfloor -1}q^{k^{2}}(-zq^{-\chi(n)})^{k}(g(\left\lfloor \frac{n}{2}\right\rfloor +k)-1)+\sum_{k=\left\lfloor n/4\right\rfloor }^{\left\lfloor (n+1)/2\right\rfloor }q^{k^{2}}(-zq^{-\chi(n)})^{k}g(\left\lfloor \frac{n}{2}\right\rfloor +k).\nonumber \end{align}
 For $\alpha>-1$ and sufficiently large $n$ with\begin{equation}
0<\frac{2q^{n/4}}{1-q}<1,\label{eq:3.71}\end{equation}
we expand \begin{align}
g(k)-1 & =\left\{ r_{2}\left(q^{\alpha+1};n\right)+1\right\} \left\{ r_{1}\left(q;n-k\right)+1\right\} \left\{ r_{1}\left(q;k\right)+1\right\} \left\{ r_{1}\left(q^{\alpha+1};n-k\right)+1\right\} -1\label{eq:3.72}\end{align}
 and estimate each term using Lemma\ref{lem:1} to get\begin{equation}
|g(\left\lfloor \frac{n}{2}\right\rfloor -k)-1|\le\frac{30q^{n/4}}{1-q},\quad0\le k\le\left\lfloor \frac{n}{4}\right\rfloor ,\label{eq:3.73}\end{equation}
 and\begin{equation}
|g(\left\lfloor \frac{n}{2}\right\rfloor +k)-1|\le\frac{30q^{n/2}}{1-q}<\frac{30q^{n/4}}{1-q},\quad k\in\mathbb{N}.\label{eq:3.74}\end{equation}
Therefore,\begin{equation}
L_{n}^{(\alpha)}(zq^{-n-\alpha};q)=\frac{(-z)^{\left\lfloor n/2\right\rfloor }\left\{ \theta_{4}(z^{-1}q^{\chi(n)};q)+e(n)\right\} }{(q;q)_{\infty}^{2}q^{\left\lfloor n/2\right\rfloor \left\lfloor (n+1)/2\right\rfloor }},\label{eq:3.75}\end{equation}
and\begin{equation}
|e(n)|\le60\theta_{3}(|z|^{-1}q^{\chi(n)};q)\left\{ \frac{q^{\left\lfloor n/4\right\rfloor ^{2}+\chi(n)\left\lfloor n/4\right\rfloor }}{|z|^{\left\lfloor n/4\right\rfloor }}+|z|^{\left\lfloor n/4\right\rfloor }q^{\left\lfloor n/4\right\rfloor ^{2}-\chi(n)\left\lfloor n/4\right\rfloor }+\frac{q^{n/4}}{1-q}\right\} \label{eq:3.76}\end{equation}
for sufficiently $n$. 

Formula \eqref{eq:3.59} and transformations for theta function imply
that

\begin{align}
 & L_{n}^{\alpha}(-\exp2\pi(u+n^{1-a}+\alpha n^{-a});\exp(-2\pi n^{-a}))\label{eq:3.77}\\
 & =\frac{\exp\left\{ \frac{\pi n^{-a}}{2}(n^{a}u+n)^{2}\right\} \left\{ 1+\mathcal{O}\left(\exp(-\pi n^{a}/2)\right)\right\} }{\sqrt{2n^{a}}\exp\pi\left\{ \pi n^{-a}/6-\pi n^{a}/6\right\} },\nonumber \end{align}
and

\begin{align}
 & L_{n}^{\alpha}(\exp2\pi(u+n^{1-a}+\alpha n^{-a});\exp(-2\pi n^{-a}))\label{eq:3.78}\\
 & =\sqrt{\frac{2}{n^{a}}}\frac{\exp\left\{ \frac{\pi n^{-a}}{2}\left(n^{a}u+n\right)^{2}\right\} \left\{ \cos\frac{\pi}{2}(n^{a}u+n)+\mathcal{O}\left(\exp(-\pi n^{a})\right)\right\} }{\exp\left\{ \pi n^{a}/24+\pi n^{-a}/6\right\} }\nonumber \end{align}
 as $n\to\infty$.

Since\begin{align}
\frac{q^{n}(q;q)_{n}}{(q^{\alpha+1};q)_{n}}w_{q\ell}(zq^{-n-\alpha}) & =-\frac{\sin(\pi\alpha)(q;q)_{\infty}}{\pi(q^{-\alpha};q)_{\infty}}\frac{z^{\alpha}q^{(1-\alpha)n-\alpha^{2}}(q;q)_{n}}{(q^{\alpha+1};q)_{n}(-zq^{-n-\alpha};q)_{\infty}}\label{eq:3.79}\\
 & =-\frac{\sin(\pi\alpha)(q;q)_{\infty}^{4}}{\pi(1-q^{-\alpha})(q,q^{1-\alpha},q^{\alpha+1};q)_{\infty}}\nonumber \\
 & \times\frac{z^{-n+\alpha}q^{n(n+3)/2-\alpha^{2}}}{(1+zq^{-\alpha})(q,-zq^{1-\alpha}-z^{-1}q^{1+\alpha};q)_{\infty}}\nonumber \\
 & \times\left\{ 1+\mathcal{O}(n^{a}\exp(-2\pi n^{1-a})\right\} \nonumber \\
 & =\frac{\exp\pi(u(2\alpha-2n-1)+2\alpha^{2}n^{-a}-n^{1-a}(n+3)-2\alpha n^{-a}-n^{-a}/6-n^{a}/3)}{\theta_{1}(\alpha n^{-a}i\vert n^{-a}i)\theta_{2}(ui+\alpha n^{-a}i\vert n^{-a}i)}\nonumber \\
 & \times\frac{in^{2a}\sin\pi\alpha}{\pi}\left\{ 1+\mathcal{O}(\exp(-2\pi n^{a})\right\} \nonumber \\
 & =\frac{n^{a}\exp(-\pi n^{-a}(n^{a}u+n)^{2}-\pi u)\left\{ 1+\mathcal{O}(\exp(-\pi n^{a})\right\} }{2\pi\exp\pi(3n^{1-a}+2\alpha n^{-a}-\alpha^{2}n^{-a}+\alpha^{2}n^{a}+n^{-a}/6+n^{a}/12)}\nonumber \end{align}
 for $n$ sufficiently large. Therefore,\begin{align}
 & \ell_{n}^{(\alpha)}(\exp2\pi(u+n^{1-a}+\alpha n^{-a});\exp(-2\pi n^{-a}))\label{eq:3.80}\\
 & =\frac{1}{\sqrt{\pi}}\frac{\exp(-\pi u/2)\left\{ \cos\frac{\pi}{2}(n^{a}u+n)+\mathcal{O}\left(\exp(-\pi n^{a})\right)\right\} }{\exp\frac{\pi}{2}(3n^{1-a}+\alpha^{2}n^{a}+n^{a}/6+2\alpha n^{-a}+n^{-a}/2-\alpha^{2}n^{-a})}\nonumber \end{align}
 as $n\to\infty$. 
\end{proof}


\begin{thebibliography}{10}
\bibitem{Akhiezer}N. I. Akhiezer, \emph{The Classical Moment Problem
and Some Related Questions in Analysis}, English translation, Oliver
and Boyed, Edinburgh, 1965.

\bibitem{Andrews1} G. E. Andrews, $q$-\textit{series: Their development
and application in analysis, number theory, combinatorics, physics,
and computer algebra}, CBMS Regional Conference Series, number 66,
American Mathematical Society, Providence, R.I.~1986.

\bibitem{Andrews2} G. E. Andrews, Ramanujan's {}``Lost\char`\"{}
Note book VIII: The entire Rogers-Ramanujan function, Advances in
Math. \textbf{191} (2005), 393--407.

\bibitem{Andrews3} G. E. Andrews, Ramanujan's {}``Lost\char`\"{}
Note book IX: The entire Rogers-Ramanujan function, Advances in Math.
\textbf{191} (2005), 408--422.

\bibitem{Andrews4} G. E. Andrews, R. A. Askey, and R. Roy, \textit{Special
Functions}, Cambridge University Press, Cambridge, 1999.

\bibitem{Deift1}P. Deift, \textit{Orthogonal Polynomials and Random
Matrices: a Riemann-Hilbert Approach}, American Mathematical Society,
Providence, 2000.

\bibitem{Deift2}P. Deift, T. Kriecherbauer, K. T-R. McLaughlin, S.
Venakides, and X. Zhou, Strong asymptotics of orthogonal polynomials
with respect to exponential weights, Comm. Pure Appl. Math. \textbf{52}
(1999), 1491--1552.

\bibitem{Gasper}G. Gasper and M. Rahman, \textit{Basic Hypergeometric
Series}, second edition Cambridge University Press, Cambridge, 2004.

\bibitem{Hayman}W. K. Hayman, On the zeros of a q-Bessel function,
Contemporary Mathematics, volume 382, American Mathematical Society,
Providence, 2005, 205--216.

\bibitem{Ismail1}M. E. H. Ismail, Asymptotics of $q$-orthogonal
polynomials and a $q$-Airy function, Internat. Math. Res. Notices
2005 No 18 (2005), 1063--1088.

\bibitem{Ismail2}M. E. H. Ismail, \textit{Classical and Quantum Orthogonal
Polynomials in one Variable}, Cambridge University Press, Cambridge,
2005.

\bibitem{Ismail3}M. E. H. Ismail and X. Li, \emph{Bounds for extreme
zeros of orthogonal polynomials}, Proc. Amer. Math. Soc. \textbf{115}
(1992), 131--140.

\bibitem{Ismail4}M. E. H. Ismail and D. R. Masson, $q$-Hermite polynomials,
biorthogonal rational functions, Trans. Amer. Math. Soc. \textbf{346}
(1994), 63--116.

\bibitem{Ismail5} M. E. H. Ismail and C. Zhang, Zeros of entire functions
and a problem of Ramanujan, Advances in Math., (2007), to appear.

\bibitem{Ismail6}M. E. H. Ismail and R. Zhang, \emph{Scaled Asymptotics
for q-Polynomials}, Comptes Rendus, Vol. 344, Issue 2, 15 January
2007, Pages 71-75.

\bibitem{Ismail7}M. E. H. Ismail and R. Zhang, \emph{Chaotic and
Periodic Asymptotics for q-Orthogonal Polynomials}, joint with Mourad
E.H. Ismail, International Mathematics Research Notices, Accepted.

\bibitem{Kajiwara}K. Kajiwara, T. Masuda, M. Noumi, Y. Ohta, Y. Yamada,
Hypergeometric solutions to the $q$-Painlev\'{e} equations, Internat.
Math. Res. Notices 47 (2004), 2497--2521.

\bibitem{Koekoek}R. Koekoek and R. Swarttouw, \emph{The Askey-scheme
of hypergeometric orthogonal polynomials and its $q$-analogues},
Reports of the Faculty of Technical Mathematics and Informatics no.
98-17, Delft University of Technology, Delft, 1998.

\bibitem{Mehta} M. L. Mehta, \textit{Random Matrices}, third edition,
Elsevier, Amsterdam, 2004.

\bibitem{Qiu}W.-Y. Qiu and R.Wong, Uniform asymptotic formula for
orthogonal polynomials with exponential weight, SIAM J. Math. Anal.\textbf{31}
(2000), 992--1029.

\bibitem{Ramanujan}S. Ramanujan, \textit{The Lost Notebook and Other
Unpublished Papers} (Introduction by G. E. Andrews), Narosa, New Delhi,
1988.

\bibitem{Rademarcher}Hans Rademarcher, \emph{Topics in Analytic Number
Theory}, Die Grundlehren der mathematischen Wissenschaften, Bd. 169,
Springer-Verlag, New York-Heidelberg-Berlin, 1973. Z. 253. 10002.

\bibitem{Szego}G. Szeg\H{o}, \textit{Orthogonal Polynomials}, Fourth
Edition, Amer. Math. Soc., Providence, 1975.

\bibitem{Wang}Z. Wang and R. Wong, Uniform asymptotics for the Stieltjes-Wigert
polynomials: the Riemann-Hilbert approach, to appear.

\bibitem{Wong} R. Wong, \textit{Asymptotic Approximations of Integrals},
Academic Press, Boston, 1989.

\bibitem{Whittaker}E. T. Whittaker and G. N. Watson, \textit{A Course
of Modern Analysis}, fourth edition, Cambridge University Press, Cambridge,
1927. 
\end{thebibliography}
\end{document}